\documentclass[nohyperref]{article}

\usepackage{pifont}
\newcommand{\xmark}{\ding{55}}

\newcommand{\defeq}{\vcentcolon=}

\newcommand{\be}{\begin{equation}}
\newcommand{\ee}{\end{equation}}
\newcommand{\bee}{\begin{eqnarray}}
\newcommand{\eee}{\end{eqnarray}}
\newcommand{\bse}{\begin{subequations}}
\newcommand{\ese}{\end{subequations}}
\newcommand\norm[1]{\left\lVert#1\right\rVert}
\newcommand\algname[1]{\textsf{#1}}
\newcommand{\Hil}{\mathcal{H}}
\newcommand{\OurAlg}{CONFIG}
\newcommand{\betac}{\beta}
\newcommand{\rev}[1]{{#1}}

\ifodd 0
    \newcommand\revReb[1]{{\color{blue}#1}}
    \newcommand{\comReb}[1]{\textbf{\color{red} (COMMENT: #1)}} 
\else
    \newcommand\revReb[1]{{#1}}
    \newcommand{\comReb}[1]{}
\fi

\newcommand{\cnj}[1]{\textcolor{red}{}} 

\usepackage{microtype}
\usepackage{graphicx}
\usepackage{subfigure}
\usepackage{booktabs} 
\usepackage{adjustbox}
\usepackage{hyperref}

\usepackage[accepted]{icml2023}

\usepackage{amsmath}
\usepackage{amssymb}
\usepackage{mathtools}
\usepackage{amsthm}

\usepackage[capitalize,noabbrev]{cleveref}

\theoremstyle{plain}
\newtheorem{theorem}{Theorem}[section]

\newtheorem{lemma}[theorem]{Lemma}
\newtheorem{corollary}[theorem]{Corollary}
\theoremstyle{definition}
\newtheorem{definition}[theorem]{Definition}
\newtheorem{assumption}[theorem]{Assumption}
\newtheorem{remark}{Remark}
\theoremstyle{remark}

\usepackage[textsize=tiny]{todonotes}

\icmltitlerunning{Constrained Efficient Global Optimization of Expensive Black-box Functions}

\begin{document}

\twocolumn[
\icmltitle{Constrained Efficient Global Optimization of Expensive Black-box Functions}




\begin{icmlauthorlist}
\icmlauthor{Wenjie Xu}{epfl,empa}
\icmlauthor{Yuning Jiang}{epfl}
\icmlauthor{Bratislav Svetozarevic}{empa}
\icmlauthor{Colin Jones}{epfl}
\end{icmlauthorlist}

\icmlaffiliation{epfl}{Automatic Control Laboratory, EPFL, Lausanne, Switzerland}
\icmlaffiliation{empa}{Urban Energy Systems Laboratory, Empa, Zurich, Switzerland}

\icmlcorrespondingauthor{Yuning Jiang}{yuning.jiang@epfl.ch}

\icmlkeywords{Machine Learning, ICML}

\vskip 0.3in
]



\printAffiliationsAndNotice{}  

\begin{abstract}
We study the problem of constrained efficient global optimization, where both the objective and constraints are expensive black-box functions that can be learned with Gaussian processes. We propose \algname{\OurAlg}~(\textbf{CON}strained ef\textbf{FI}cient \textbf{G}lobal Optimization), a simple and effective algorithm to solve it. Under certain regularity assumptions, we show that our algorithm enjoys the same cumulative regret bound as that in the unconstrained case and similar cumulative constraint violation upper bounds. For commonly used M\'atern and Squared Exponential kernels, our bounds are sublinear and allow us to derive a convergence rate to the optimal solution of the original constrained problem. In addition, our method naturally provides a scheme to declare infeasibility when the original black-box optimization problem is infeasible. Numerical experiments on sampled instances from the Gaussian process, artificial numerical problems, and a black-box building controller tuning problem all demonstrate the competitive performance of our algorithm. \revReb{Compared to the other state-of-the-art methods, our algorithm significantly improves the theoretical guarantees, while achieving competitive empirical performance.}  
\end{abstract}

\section{Introduction}
Global optimization of expensive black-box functions is a pervasive problem in science and engineering. In black-box optimization problems, we typically need to sequentially evaluate different candidate solutions to find the optimal one without explicit gradient information. For example, hyperparameter tuning for machine learning models~\cite{bergstra2012random,snoek2015scalable}, control system performance optimization~\cite{bansal2017goal,xu2022vabo} and drug discovery or materials design~\cite{negoescu2011knowledge,frazier2016bayesian} can all be formulated as black-box optimization problems.  

\rev{In these applications, we typically know little about the shape of the black-box functions, which can be non-convex and multi-modal. As a result, a global optimization method is needed.
In many applications, e.g., tuning the hyperparameters of a deep neural network, each function evaluation can take a large number of resources such as energy, time, or GPU computation, so we want to use as few function evaluations as possible. Standard global optimization methods, such as genetic algorithms or simulated annealing~\cite{pardalos2013handbook}, typically require a large number of samples and thus can not be applied to the optimization of expensive black-box functions.   
} 
\rev{Therefore, Gaussian process based efficient global optimization\footnote{Also known as Gaussian process optimization~(e.g., in~\cite{srinivas2012information}), Bayesian optimization~(e.g., in~\cite{frazier2018tutorial,hernandez2016general}) or kriging~(e.g., in~\cite{jeong2005efficient}). }\cite{jones1998efficient,shahriari2015taking}, as a sample-efficient method to solve black-box optimization problems~\cite{xu2022lower}, has recently been becoming more and more popular. The general idea of efficient global optimization is to build a surrogate model, typically by Gaussian process regression, of the black box function from samples and then to choose a search direction using this model. This problem, when formulated as a (kernelized) multi-armed bandit problem, is widely studied in the  literatures~\cite{audibert2010best,srinivas2012information,soare2014best,chowdhury2017kernelized,amani2019linear, zhou2022kernelized}. 
}

A key challenge in many applications is black-box constraints. For example, when tuning the temperature controller parameters of a building, we need to minimize the energy consumption while keeping the occupants' comfort at user-defined levels, where both the energy and the comfort are black-box functions of the controller parameters. 
To handle the unknown black-box constraints, a variety of efficient global optimization methods with constraints have been developed. We can roughly classify them into different groups based on whether constraint violations are allowed during the optimization process. For the setting where no constraint violations are tolerated, a group of safe optimization methods has been developed~\cite{sui2015safe,wu2016conservative,sui2018stagewise,turchetta2020safe,moradipari2020stage, amani2020regret}. However, these safe algorithms require a set of feasible solutions known a priori, while in practice, we may not get access to such feasible solutions or even do not know the feasibility of the problem. Furthermore, due to these hard safety-critical constraints, these algorithms may become stuck at a local minimum. The other extreme is where constraint violations are allowed to occur during the optimization process without penalty, for which there exists generic constrained Bayesian~(efficient global) optimization methods~\cite{schonlau1997computer,sasena2002exploration,basudhar2012constrained,bagheri2017constraint, zhang2021two}. Along this line of research, a group of popular methods~\cite{gardner2014bayesian,gelbart2014bayesian} encodes the constraint information into the acquisition function~(e.g., Constrained Expected Improvement). However, there are no theoretical guarantees on optimality, constraint violation incurred during optimization, or convergence rates for this class of methods. A recent work of~\cite{inatsu2022bayesian} considers Bayesian optimization for distributionally robust chance-constrained problems but can not provide cumulative regret/violation bounds.

More recently, a third group of approaches has been studied~(e.g.,~\cite{xu2022vabo,zhou2022kernelized}), where constraint violations are allowed, but limits are placed on the total violation incurred during the optimization process. In this setting, two more recent works adopt a penalty-function approach~\cite{lu2022no} and a primal-dual approach~\cite{zhou2022kernelized} to solve constrained efficient global optimization problems. The common ideas behind these two approaches are the addition of a penalizing term of the constraint violation in the objective and the transformation of the constrained problem to an unconstrained one. However, both of these methods require choosing values for some critical parameters~(e.g., penalty coefficient~\cite{lu2022no} and dual update step size~\cite{zhou2022kernelized}). The performance of these methods relies on the chosen parameters heavily, and hence these methods can take significant parameter tuning effort for implementation. Additionally, the cumulative violations considered in~\cite{zhou2022kernelized} are the violations of the cumulative constraint values, not the real total violations. Such a form of cumulative violation bounds is weak and can not rule out the case where points with severe violations and small constraint values are alternatingly sampled while keeping a low cumulative constraint value. The work of~\cite{lu2022no} only provides a penalty-based regret bound, which is defined as the regret plus a penalty weight times violations, which does not directly lead to separate bounds on regret and violation since the penalty weight is unknown beforehand. We also notice that a concurrent work~\cite{guo2022rectified} extends the fixed penalty to an increasing penalty sequence. However, it requires a delicate design of the penalty sequence and may result in numerical issues when the penalty is large. Last but not least, all the existing methods~\cite{sui2015safe,gelbart2014bayesian,zhou2022kernelized,lu2022no}, either require an initial feasible solution to be given or lack the scheme to detect and report infeasibility when the original problem is infeasible. This is key, as, in practice, black-box optimization problems may be infeasible, and a scheme to detect and report infeasibility is desired.

\begin{figure*}
    \centering
    \includegraphics[width=0.8\textwidth]{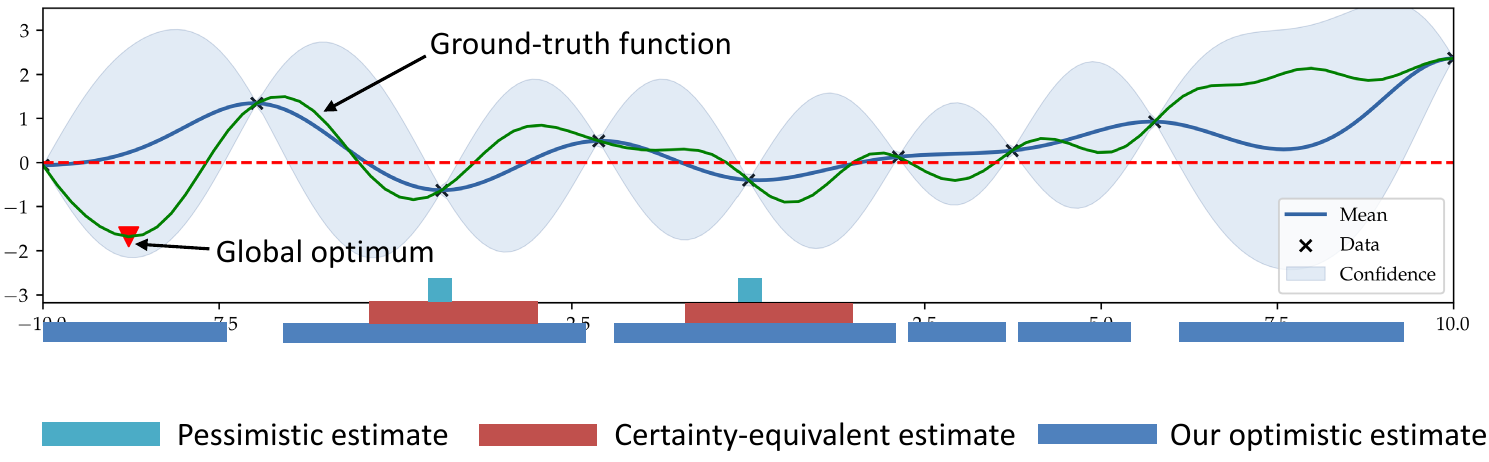}
    \caption{\revReb{A simple demonstration of different ways to construct the estimate of the feasible set, where we minimize an unknown function subject to that it is non-positive. Pessimistic estimate can guarantee feasibility with high probability. However, it can be easily seen that the pessimistic estimate is very tiny and can easily be empty if the initial several points are infeasible. With the certainty-equivalent estimate, the algorithm may easily get stuck due to lack of exploration. To the contrast, we construct an optimistic estimate of the feasible set using lower confidence bound functions, which contains the ground-truth feasible set with high probability. }}
    \label{fig:config_demo}
\end{figure*}
\begin{table*}[htbp!]
\small
\renewcommand{\arraystretch}{1.2}
\centering

{\caption{The comparison of our method to existing state-of-the-art constrained efficient global optimization~(or kernelized multi-armed bandits with constraints) methods.  } 
\label{tab:comparison}
\begin{tabular}{|c|c|c|c|c|}
\hline 
\textbf{Works} & \begin{tabular}{@{}c@{}}\textbf{Cumulative}\\ \textbf{Regret Bound}\end{tabular}&\begin{tabular}{@{}c@{}}\textbf{Cumulative}\\ \textbf{Violation Bound}\end{tabular} & \begin{tabular}{@{}c@{}}\textbf{Optimality}\\{\textbf{Guarantee}}\end{tabular} & \begin{tabular}{@{}c@{}}\textbf{Infeasibility}\\\textbf{Declaration}\end{tabular}\\
\hline 
\begin{tabular}{@{}c@{}}
\textsf{SafeOPT}~
\cite{sui2015safe}, etc.
\end{tabular}& \xmark & No violation & Local convergence& \xmark \\
\hline 
\begin{tabular}{@{}c@{}}
\textsf{Constrained EI}~
\cite{gardner2014bayesian}, etc.
\end{tabular}& \xmark& \xmark & \xmark &\xmark \\
\hline
\begin{tabular}{@{}c@{}}
\textsf{Primal-Dual}
~
\cite{zhou2022kernelized}, etc.
\end{tabular} & \checkmark & weak-form& \xmark & \xmark \\
\hline 
\begin{tabular}{@{}c@{}}
\textsf{Penalty function method}\\
\cite{lu2022no}, etc.
\end{tabular}& \multicolumn{2}{c|}{Penalty-based Regret Bound} & \begin{tabular}{@{}c@{}}Global convergence\\ with penalty-dependent rate\end{tabular} & \xmark \\
\hline 
\begin{tabular}{@{}c@{}}
\textsf{CONFIG}~
(This paper)
\end{tabular} & \checkmark& \checkmark & \begin{tabular}{@{}c@{}}Global convergence\end{tabular}  & \checkmark \\
\hline 
\end{tabular}
}  
\end{table*}

\textbf{Our contributions.} We propose \textsf{\OurAlg}, a simple and effective algorithm for the constrained efficient global optimization of expensive black-box functions by exploiting the principle of \emph{optimism in the face of uncertainty}. Specifically, in each step, our algorithm solves an auxiliary constrained optimization problem with lower confidence bound~(LCB) surrogates as both objective and constraints in order to choose the next sample. \revReb{Fig.~\ref{fig:config_demo} demonstrates the idea of the optimism in the feasible set. Intuitively, our algorithm strategically trades potential violation for faster learning of both the objective function and the constraints. }
%
Under certain regularity assumptions, we show that the algorithm enjoys the same cumulative regret bound as that in the unconstrained case. In addition, we show that the cumulative constraint violation has an upper bound in terms of maximum information gain, which is similar to the cumulative regret bounds of the objective function. For the commonly used M\'atern and Squared Exponential kernels, our bounds are sublinear and allow us to derive a convergence rate to the optimal solution of the original constrained problem. Furthermore, our method naturally provides a scheme to declare infeasibility when the original black-box optimization problem is infeasible. \revReb{To the sharp contrast, existing algorithms for constrained Bayesian optimization, including the popular CEI method~\cite{gardner2014bayesian} and the recently introduced primal-dual method~\cite{zhou2022kernelized,xu2023primal}, can only provide no or weak guarantees on the cumulative regret/violation bounds or the convergence to the optimal solution (See Tab.~\ref{tab:comparison} for more details). Additionally, none of the existing methods can detect infeasibility when the original problem is infeasible.} We then run numerical experiments to corroborate the effectiveness of our algorithm. Numerical results show that our algorithm achieves performance competitive with the state-of-the-art constrained expected improvement~(CEI) algorithm in terms of both optimality and constraint violations. \revReb{As such, our algorithm significantly improves the theoretical guarantees, while achieving competitive empirical performance, as compared to the existing state-of-the-art constrained efficient global/Bayesian optimization methods. } 

\section{Notation and Preliminaries}
This section introduces the class of constrained global optimization problems studied in this paper and recaps existing results from the field of efficient global optimization.  
\subsection{Problem statement}
We consider the following black-box optimization problem, 
\begin{subequations}
\label{eqn:prob_form}
\begin{eqnarray}
\min_{x\in\mathcal{X}}\quad   && f(x)\enspace,\label{eqn:obj} \\
\text{subject to}\quad && g_i(x)\leq 0\enspace,\quad \forall i\in[N],\label{eqn:constraint}
\end{eqnarray}
\end{subequations}
where $\mathcal{X}\subset\mathbb{R}^d$ is a known set of candidate solutions~(e.g., a hyperbox) with dimension $d$ and $[N]:=\{1,2,\dots, N\}\subset\mathbb Z$ denotes the set of integers from $1$ to $N$. Both $f$ and $g_i, i\in[N]$ are unknown black-box functions.

\begin{remark}
In Prob.~\ref{eqn:prob_form}, we only consider inequality constraints because an equality constraint $g(x)=0$ can be equivalently transformed into two inequality constraints $g(x)\leq0$, $-g(x)\leq0$.
\end{remark}
To proceed, we make the following commonly used assumptions in the field of efficient global optimization. 
\begin{assumption}
Set $\mathcal{X}$ is compact. 
\label{assump:support_set}
\end{assumption}
Assumption~\ref{assump:support_set} is naturally satisfied in a large set of applications. For example, when tuning the hyperparameters of a machine learning model, we can restrict the parameters to lie inside a hyperbox based on prior knowledge.  
\begin{assumption}
\label{assump:bounded_norm}
$f\in\Hil_{k_0}$, $g_i\in\Hil_{k_i}, i\in[N]$, where $k_i:\mathbb{R}^d\times\mathbb{R}^d\to \mathbb{R}, i\in\{0\}\cup[N]$ are symmetric, positive-semidefinite kernel functions and $\Hil_{k_i}$ are their corresponding reproducing kernel Hilbert spaces~(RKHSs, see~\cite{scholkopf2001generalized}). Furthermore, we assume $\|f\|_{k_0}\leq B_0$ and $\|g_i\|_{k_i}\leq B_i, i\in[N]$, where $\|\cdot\|_{k_i}$ is the norm induced by the inner product in the corresponding reproducing kernel Hilbert space. 
\end{assumption}
Assumption~\ref{assump:bounded_norm} requires that the underlying functions are regular in the sense of having a bounded norm in an RKHS. In existing literature of efficient global optimization~(also known as Bayesian optimization or Gaussian process optimization), having a bounded norm in an RKHS is a common assumption~(e.g.,~\cite{srinivas2012information,zhou2022kernelized}). 

\begin{assumption}
\label{assump:feasibility}
Prob.~\eqref{eqn:prob_form} is feasible and has an optimal solution $x^*$ with the optimal value $f^*=f(x^*)$. 
\end{assumption}
With Assumption~\ref{assump:feasibility}, we focus on the case that the original black-box optimization problem is feasible. We will separately discuss how our algorithm can handle infeasible problem instances later in Sec.~\ref{sec:inf_dec}. 

We consider a zero-order feedback setting, where the algorithm can sequentially query noisy evaluations of different candidate solutions to learn the optimal feasible solution. At step $t$, with query point $x_t$, we get noisy evaluations,
\bse
\begin{align}
y_{0,t}&=f(x_t)+\xi_{0,t}\enspace,\\
y_{i,t}&=g_i(x_t)+\xi_{i,t}\enspace,\quad i\in[N],
\end{align}
\ese
where $\xi_{i,t},i\in\{0\}\cup[N], t\geq1$ are i.i.d $\sigma$-sub-Gaussian noise with fixed $\sigma>0$. 
We also use $X_t\defeq(x_1, x_2, \cdots, x_t)$ to denote the sequence of sampled points up to step $t$.

\subsection{Gaussian process surrogates}
To design our algorithm, we use a Gaussian process surrogate model to learn the unknown functions. As in~\cite{chowdhury2017kernelized}, we artificially introduce a Gaussian process $\mathcal{GP}(0, k(\cdot, \cdot))$ for the surrogate modeling of the unknown black-box function $f$. We also adopt an \emph{i.i.d} Gaussian zero-mean noise model with noise variance $\lambda>0$.
\begin{remark}
Note that this Gaussian process model is \emph{only} used to derive posterior mean functions, covariance functions, and maximum information gain for the purpose of algorithm design. We introduce the artificial Gaussian process only to make the algorithm more interpretable. It does not change our set-up that $f$ is a deterministic function and that the noise is only assumed to be sub-Gaussian. Indeed, we are considering the agnostic setting introduced in~\cite{srinivas2012information}.
\end{remark}
We introduce the following functions for $x, x^{\prime}$,
\bse
\label{eq:mean_cov}
\begin{align}
\mu_{0,t}(x) &=k_{0}(x_{1:t}, x)^\top\left(K_{0,t}+\lambda I\right)^{-1} y_{0, 1: t}\enspace, \\
k_{0,t}\left(x, x^{\prime}\right) &=k_{0}\left(x, x^{\prime}\right)\\\notag
&-k_{0}(x_{1:t}, x)^\top\left(K_{0,t}+ \lambda I\right)^{-1} k_{0}\left(x_{1:t}, x^{\prime}\right), \\
\sigma_{0,t}^2(x) &=k_{0,t}(x, x)\enspace,
\end{align}
\ese
where $k_{0}(x_{1:t}, x)=[k_0(x_1, x), k_0(x_2, x),\cdots, k_0(x_t, x)]^\top$, $K_{0,t}=(k_0(x,x^\prime))_{x,x^\prime\in X_t}$ and $y_{0,1:t}=[y_{0,1}, y_{0, 2},\cdots,y_{0,t}]^\top$. Similarly, we can derive $\mu_{i,t}(\cdot), k_{i,t}(\cdot, \cdot), \sigma_{i,t}(\cdot)$, $\forall i\in[N]$ for the constraints. We also introduce the maximum information gain for the objective $f$ as in~\cite{srinivas2012information},
\bee
\label{eq:max_inf_gain}
\gamma_{0,t}:=\max _{X \subset \mathcal{X} ;|X|=t} \frac{1}{2} \log \left|I+\lambda^{-1}K_{0,X}\right|,
\eee
where $K_{0,X}=(k_0(x,x^\prime))_{x,x^\prime\in X}$. 
Similarly, we can introduce the maximum information gain $\gamma_{i,t},\forall i\in[N]$ for the constraints. 

\begin{lemma}[Theorem 2,~\cite{chowdhury2017kernelized}]
Let Assumptions~\ref{assump:support_set} and \ref{assump:bounded_norm} hold. For \rev{any $\delta\in(0, 1)$}, with probability at least $1-\delta/(N+1)$, the following holds for all $x \in \mathcal{X}$ and $1\leq t \leq T$, $T\in\mathbb{N}$, 
\begin{align}
&\left|\mu_{0,t-1}(x)-f(x)\right|\\\notag \leq&\left(B_0+\sigma \sqrt{2\left(\gamma_{0,t-1}+1+\ln ((N+1) / \delta)\right)}\right) \sigma_{0,t-1}(x),
\end{align}
where $\mu_{0,t-1}(x), \sigma^2_{0,t-1}(x)$ and $\gamma_{0,t-1}$ are as given in Eq.~\eqref{eq:mean_cov} and Eq.~\eqref{eq:max_inf_gain}, and $\lambda$ set to be $1+\frac{2}{T}$.~($\mu_{0,0}$ and $\sigma_{0,0}$ are the prior mean and standard deviation, and $\gamma_{0,0}=0$.) \label{lem:conf_int} 
\end{lemma}
\noindent
\rev{We remark that we replace the `$\delta$' in Theorem 2 of~\cite{chowdhury2017kernelized} by $\frac{\delta}{N+1}$, which will be useful to derive the confidence interval with $1-\delta$ probability in Corollary~\ref{cor:hpb}. 
}
We can derive similar confidence intervals for the constraint functions. To facilitate the following algorithm design and discussion, we introduce the lower confidence and upper confidence bound functions,
\begin{definition}
\label{def:lw_ub}
For $i\in\{0\}\cup[N]$, $x\in\mathcal{X}$ and $t\in[T]$, 
\bse
\begin{align}
l_{i,t}(x)&\defeq\mu_{i,t-1}(x)-\betac_{i,t}\sigma_{i,t-1}(x)\enspace,\\
u_{i,t}(x)&\defeq\mu_{i,t-1}(x)+\betac_{i,t}\sigma_{i,t-1}(x)\enspace,
\end{align}
\ese
where $\betac_{i,t}=B_i+\sigma \sqrt{2\left(\gamma_{i,t-1}+1+\ln ((N+1) / \delta)\right)}$.
\label{def:hpb_int}
\end{definition}
We then have the following corollary, 
\begin{corollary}
\label{cor:hpb}
Let Assumptions~\ref{assump:support_set} and \ref{assump:bounded_norm} hold. 
With probability at least $1-\delta$, the following holds for all $x\in\mathcal{X}$ and $1\leq t\leq T$,
\begin{align}
g_i(x)&\in[l_{i,t}(x),u_{i,t}(x)]\enspace,\quad\forall i\in[N]\\ 
\textrm{ and }\quad f(x)&\in[l_{0,t}(x),u_{0,t}(x)].
\end{align}
\end{corollary}

\subsection{Performance metric}
Before we design our algorithm, we introduce the performance metric. We are interested in minimizing the gap of $f(x_t)$ to the optimal value $f^*=\min_{x\in \mathcal{X}, g_i(x)\leq0, \forall i\in[N]}f(x)$, i.e., the instantaneous regret. 
Instantaneous regret at step $t$ is defined as, 
\bee
r_t = f(x_t)-f^*,
\eee
where $x_t$ is the point queried by our algorithm in step $t$. 
Since $f^*$ is the constrained optimal value, we may sample some infeasible points with even smaller objective values than $f^*$, we further define the positive regret as,
\bee 
r_t^+=[f(x_t)-f^*]^+,
\eee
where $[\cdot]^+\defeq\max\{0, \cdot\}$. We are also interested in the constraint violation, defined as,
\bee
v_{i,t} = [g_i(x_t)]^+.
\eee
We then introduce cumulative regret that measures the extra cost incurred during the running of the algorithm, as in existing kernelized multi-armed bandit literature~(e.g.,~\cite{srinivas2012information,chowdhury2017kernelized}).  
\begin{definition}[Cumulative-Regret]
We define the cumulative-regret as
\bee
{R}_T=\sum_{t=1}^T(f(x_t)-f^*).
\eee 
To facilitate the following discussion, we further define the cumulative positive regret,
\bee
{R}_T^+=\sum_{t=1}^T[f(x_t)-f^*]^+.
\eee 
\end{definition}
We further define the cumulative violation.
\begin{definition}[Cumulative-Violation]
We define the cumulative-violation as
\bee
\mathcal{V}_{i,T}=\sum_{t=1}^T[g_i(x_t)]^+, \forall i\in[N].
\eee   
\end{definition}
Cumulative regret/violations measure the performance during the running of the algorithm. Besides them, we are also interested in the convergence speed to the constrained global optimum. Intuitively, small average regret/violations indicate the existence of a sample with small regret/violation. 

\section{Algorithm}
We propose the lower confidence bound based \textbf{CON}strained ef\textbf{FI}cient \textbf{G}lobal (\textbf{CONFIG}) optimization algorithm to solve our problem as shown in Alg.~\ref{alg:lcb2}~\footnote{Knowledge of $T$ is assumed. Standard `doubling trick' can be applied to extend to the setting without knowing $T$.}. 
\begin{algorithm}[htbp!]
\caption{Lower confidence bounds based \textbf{CON}strained ef\textbf{FI}cient \textbf{G}lobal (\textbf{CONFIG}) optimization algorithm.}
\begin{algorithmic}[1]
\normalsize
\FOR{$t\in[T]$}
\IF{$\exists i\in[N]$, such that $\min_{x\in\mathcal{X}}l_{i,t}(x)>0$}
\STATE \textbf{Declare infeasibility}. \label{alg_line:declare_inf}
\ENDIF
\STATE Find $x_t\in{\arg\min}_{x\in \mathcal{X}} l_{0,t}(x)$\\\enspace\qquad${\textbf{ subject to }} l_{i,t}(x)\leq0,\forall i\in[N]$.\label{alg_line:aux_prob}
\STATE Make noisy evaluations of $f$ and $g_i,i\in[N]$ at $x_t$.
\STATE Update Gaussian process posterior mean and covariance with the new evaluations added. 
\ENDFOR
\end{algorithmic}
\label{alg:lcb2}
\end{algorithm}
Our algorithm is conceptually simple based on the \emph{optimism in the face of uncertainty} and only requires one to solve an auxiliary problem with the black-box objective and constraints replaced by their lower confidence bound surrogates as shown in line~\ref{alg_line:aux_prob} of Alg.~\ref{alg:lcb2}. Since we consider problems with expensive-to-evaluate functions, the cost of solving the auxiliary problem, which may itself be non-convex, is much smaller compared to solving the original black-box optimization problem when the dimension of the input space is small to medium size. To solve the auxiliary problem, we can use a pure grid search method when the dimension is small~(e.g., $\leq5$). For larger dimensions, one can, for example, start from multiple different points and apply gradient-based methods.

\rev{
Note that in line~\ref{alg_line:declare_inf} of Alg.~\ref{alg:lcb2}, we may declare infeasibility if the auxiliary problem in line~\ref{alg_line:aux_prob} is infeasible. When the original problem~\eqref{eqn:prob_form} is feasible, in some very rare cases, we may find the auxiliary problem infeasible due to the potentially unbounded noise. Fortunately, we will show that with guaranteed high probability, line~\ref{alg_line:declare_inf} of Alg.~\ref{alg:lcb2} is never declared when the original problem is feasible and Assumptions~\ref{assump:support_set} and \ref{assump:bounded_norm} hold.   
}
\cnj{You've got Declare Infeasibility as a possible outcome of this algorithm, but one of the first assumptions made is that the problem is feasible}

\section{Analysis}
We now analyze the performance of Alg.~\ref{alg:lcb2}. 
All the proofs to the results are put in the appendix.
We first introduce a lemma to bound the instantaneous regret and violation,
\begin{lemma}
Let Assumptions~\ref{assump:support_set}, \ref{assump:bounded_norm}, and \ref{assump:feasibility} hold. With probability at least $1-\delta$, we have first that infeasibility is {never} declared in line~\ref{alg_line:declare_inf} of the Alg.~\ref{alg:lcb2}, and for all $1\leq t\leq T$,
\begin{align}
r_t&\leq r_t^+\leq2\betac_{0,t}\sigma_{0,t-1}(x_t)\enspace,\\
v_{i,t}&\leq2\betac_{i,t}\sigma_{i,t-1}(x_t)\enspace,\quad i\in[N].
\end{align}
\label{lem:bound_ind_rv}
\end{lemma}
We can also bound the cumulative sum of the posterior standard deviation,
\begin{lemma}[Lemma 4,~\cite{chowdhury2017kernelized_arxiv}]
If $x_1, x_2, \cdots, x_T$ are the points selected by Alg.~\ref{alg:lcb2}, then,
\bee
\label{lem:bound_cumu_sd}
\sum_{t=1}^T \sigma_{i, t-1}\left(x_t\right) \leq \sqrt{4(T+2) \gamma_{i,T}}.
\eee
\end{lemma}
The bound in Eq.~\eqref{lem:bound_cumu_sd} is in terms of maximum information gain. To derive a kernel-specific bound, we need to estimate the maximum information gain. A set of upper bounds for the maximum information gains were given by~\revReb{\cite{srinivas2012information,vakili2021information}}, which we restate in the maximum information gain column of Tab.~\ref{tab:kernel_spec_bounds}. 
\begin{table*}[htbp!]
\caption{Kernel-specific cumulative regret/violation bounds and the convergence rate to the optimal solution. See Appendix~\ref{app:ex_kerfuncs} for the specific kernel functions. The column of maximum information gain shows the results of~\revReb{\cite{srinivas2012information,vakili2021information}}. In the table, $d$ represents the dimension of the input space, and $\nu$ represents the smoothness parameter of the M\'atern kernel.} 
\label{tab:kernel_spec_bounds}
\centering
\begin{tabular}{|c|c|c|c|c|}
\hline
\textbf{Kernel} & \textbf{Maximum Information Gain}& \textbf{Cumulative Regret/Violation} & \textbf{Convergence Rate} \\
\hline 
\textbf{Linear} &$\mathcal{O}(d\log T)$ &$\mathcal{O}\left(d\log T\sqrt{T}\right)$  &$\mathcal{O}\left(\frac{(N+1)d\log T}{\sqrt{T}}\right)$  \\
\hline
\begin{tabular}{@{}c@{}}
\textbf{Squared }\\[-0.1cm]
\textbf{Exponential}
\end{tabular}
& $\mathcal{O}((\log T)^{d+1})$ &$\mathcal{O}\left((\log T)^{d+1}\sqrt{T}\right)$ &$\mathcal{O}\left(\frac{(N+1)(\log T)^{d+1}}{\sqrt{T}}\right)$  \\
\hline 
\begin{tabular}{@{}c@{}}
\textbf{M\'atern}\\
\revReb{$\left(\nu>\frac{d}{2}\right)$}
\end{tabular}
&\revReb{$\mathcal{O}(T^{\frac{d }{2 \nu+d}}\log^{\frac{2\nu}{2\nu+d}}(T))$} & \revReb{$\mathcal{O}(T^{\frac{2\nu+3d }{4 \nu+2d}}\log^{\frac{2\nu}{2\nu+d}}(T))$}  &\revReb{$\mathcal{O}\left((N+1)T^{-\frac{2\nu-d}{4 \nu+2d}}\log^{\frac{2\nu}{2\nu+d}}(T)\right)$}   \\ 
\hline
\end{tabular}
\end{table*}
We can now state our main theorem,
\begin{theorem}
\label{thm:R_bound}
Under Assumptions~\ref{assump:support_set}, \ref{assump:bounded_norm} and \ref{assump:feasibility}, we have, with probability at least $1-\delta$, the sample points of Algorithm 1 satisfy, 
\bse
\label{eqn:cumu_reg_bound}
\begin{align}
R_T&\leq R_T^+\leq 4\betac_{0,T}\sqrt{(T+2)\gamma_{0,T}}=\mathcal{O}(\gamma_{0,T}\sqrt{T})\enspace,\\
\mathcal{V}_{i,T}&\leq 4\betac_{i,T}\sqrt{(T+2)\gamma_{i,T}}=\mathcal{O}(\gamma_{i,T}\sqrt{T})\enspace,\quad\forall i\in[N].
\end{align}
\ese
\end{theorem}

Interestingly, the regret bound in Thm.~\ref{thm:R_bound} is exactly the \emph{same} as the regret bound in the unconstrained case shown in~\cite{chowdhury2017kernelized_arxiv}. Furthermore, the constraint violation bound is similar to the regret bound. Based on Thm.~\ref{thm:R_bound}, we can derive a convergence rate to the optimal feasible solution.
\begin{theorem}
Under Assumptions~\ref{assump:support_set}, \ref{assump:bounded_norm} and \ref{assump:feasibility}, we have, with probability at least $1-\delta$, that there exists $\tilde{x}_T\in\{x_1, x_2,\cdots, x_T\}$, such that, $\forall j\in[N]$, 
\begin{subequations}
\footnotesize
\begin{align}
f(\tilde{x}_T)-f^*&\leq \frac{4\sqrt{T+2}\sum_{i=0}^N\betac_{i,T}\sqrt{\gamma_{i,T}}}{T}=\mathcal{O}\left(\frac{\sum_{i=0}^N\gamma_{i,T}}{\sqrt{T}}\right)\enspace,\\
[g_j(\tilde{x}_T)]^+&\leq \frac{4\sqrt{T+2}\sum_{i=0}^N\betac_{i,T}\sqrt{\gamma_{i,T}}}{T}=\mathcal{O}\left(\frac{\sum_{i=0}^N\gamma_{i,T}}{\sqrt{T}}\right).
\end{align}
\end{subequations}
\label{thm:converge_opt}
\end{theorem}

Thm.~\ref{thm:converge_opt} gives the convergence rate to the constrained optimal value in terms of maximum information gain for general kernels. We can now apply Thm.~\ref{thm:R_bound} and Thm.~\ref{thm:converge_opt} to derive kernel-specific results. Specifically, if $k_i, i\in\{0\}\cup[N]$ share the same kernel types, we can list the kernel-specific bounds and the convergence rates in Tab.~\ref{tab:kernel_spec_bounds}.  

\section{Infeasibility Declaration}
\label{sec:inf_dec}
Thm.~\ref{thm:R_bound} and Thm.~\ref{thm:converge_opt} are derived under the assumption that the problem is feasible. However, the underlying problem may indeed sometimes be infeasible and therefore, we need an infeasibility detection scheme. In the following, we show that if the original problem is infeasible, we are able to declare infeasibility in line~\ref{alg_line:declare_inf} of Alg.~\ref{alg:lcb2} and detect the original problem's infeasibility with high probability.               

For infeasibility declaration, we try to avoid `false positives', where a feasible problem is declared as infeasible, and `false negatives', where an infeasible problem is not declared to be infeasible. 
Lem.~\ref{lem:bound_ind_rv} already shows that if the original black-box optimization problem is feasible, then with high probability infeasibility is not declared. The following theorem complements Lem.~\ref{lem:bound_ind_rv} by showing that if the original problem is infeasible, our algorithm is guaranteed to declare infeasibility within finite steps with high probability.  
\begin{theorem}
Let Assumptions~\ref{assump:support_set} and \ref{assump:bounded_norm} hold. Assume that Prob.~\eqref{eqn:prob_form} is infeasible, that is, 
\[
\min_{x\in\mathcal{X}}\max_{i\in[N]}\;\;g_i(x)=\epsilon>0,
\] 
and $\lim_{T\to\infty}\frac{\gamma_{i,T}}{\sqrt{T}}=0,\forall i\in[N]$. Given a desired confidence level $\delta\in(0,1)$, Alg.~\ref{alg:lcb2} 
will declare infeasibility on line~\ref{alg_line:declare_inf} within a number of steps equal to 
\bee
\overline{T}=\min_{T\in\mathbb{N}_+}\left\{T\left|\sum_{i\in[N]}\frac{\gamma_{i,T}}{\sqrt{T}}< C\epsilon\right\}\right., 
\eee
with probability at least $1-\delta$, where $C$ is a constant independent of $T$. 
\label{thm:declare_inf}
\end{theorem}
\revReb{
\begin{remark}
Our algorithm design and analysis in this paper focus on the case where $T$ is given. In cases where $T$ is not given or not large enough to declare the infeasibility, we can apply the `doubling trick'. The basic idea is that we can start running a round of steps with $T=1$ and double $T$ after each round. As $T$ grows larger and larger, it will finally be larger than $\overline{T}$ and enough for infeasibility declaration.  
\end{remark}
}
\revReb{
\begin{remark}
Note that in Thm.~\ref{thm:declare_inf}, we assume that $\lim_{T\to\infty}\frac{\gamma_{i,T}}{\sqrt{T}}=0$. From Tab.~\ref{tab:kernel_spec_bounds}, we can see that for linear and squared exponential kernels, it is satisfied and for M\'atern kernel, we require that $\nu>\frac{d}{2}$.
\end{remark}
}
\section{Experiments}
To demonstrate the effectiveness of \algname{\OurAlg}, we run experiments over sampled instances from the Gaussian process, artificial numerical instances, and a room temperature controller tuning problem. We compare our method to the state-of-the-art Gaussian process optimization baseline methods with constraints consideration. The baselines include \textsf{SafeOPT}~\cite{sui2015safe}, \textsf{CEI}~(\textbf{C}onstrained \textbf{E}xpected \textbf{I}mprovement) based method~\cite{gelbart2014bayesian,gardner2014bayesian}, primal-dual algorithm~\cite{zhou2022kernelized} and EPBO~(exact penalty Bayesian optimization)~\cite{lu2022no}. The experiments are implemented in \textsf{python}, based on the package \textsf{GPy}~\cite{gpy2014}. 

\revReb{\textbf{Choice of Parameters.} Definition~\ref{def:lw_ub} gives a rigorous way of selecting the coefficient $\beta_{i,t}$. Since $B_i$ and $\sigma$ may be hard to obtain in practice, $\beta_{i,t}$ can usually be set as a constant. The choice of $\beta_{i,t}$ is a key factor that impacts the performance of different algorithms. Empirically, if $\beta_{i,t}$ is set too small, exploration may be insufficient and we may miss the global optimal solution. If $\beta_{i,t}$ is set too large, the algorithm overly explores, converges slower and may suffer from high cumulative regret/violation before converging to the optimal solution. In our experiments shown in this section, manually setting $\beta_{i,t}=3$ works well. We also set $\lambda=0.05^2$ as the noise variance for the Gaussian process modelling. We use the common squared exponential kernel functions. For Sec.~\ref{subsec:art_nume} and Sec.~\ref{subsec:tune_P}, we randomly sample a few points and maximize the likelihood function to get the hyperparameters of the kernel.}  

\revReb{\textbf{Computational Time.} In our experiments, all the problems have low-dimensional input~($\leq3$). So we use pure grid search to solve the auxiliary problem for different algorithms. Therefore, the computational time is almost the same for different algorithms. With the explicit expression of the lower confidence bounds, the computation cost of solving the auxiliary problems is much cheaper than evaluating the expensive functions, e.g., energy consumption and discomfort of a simulated building in Sec. 6.3.    
}  

\subsection{Numerical Results over Sampled Instances from Gaussian Process}

In this section, we consider the constrained problem with both the objective and the constraint sampled from a Gaussian process. We use the squared exponential kernel as shown in~\eqref{eq:se}.
\bee
k(x,y) = \sigma_\mathrm{SE}^2\exp{\left\{-\frac{\norm{x-y}^2}{l^2}\right\}}, \label{eq:se} 
\eee
where $\sigma_{\mathrm{SE}}^2=2.0$ and $l=1.0$. 
We remove those instances with empty feasible sets. We sampled $48$ instances in total. 


\begin{figure}[htbp!]
\centering    \includegraphics{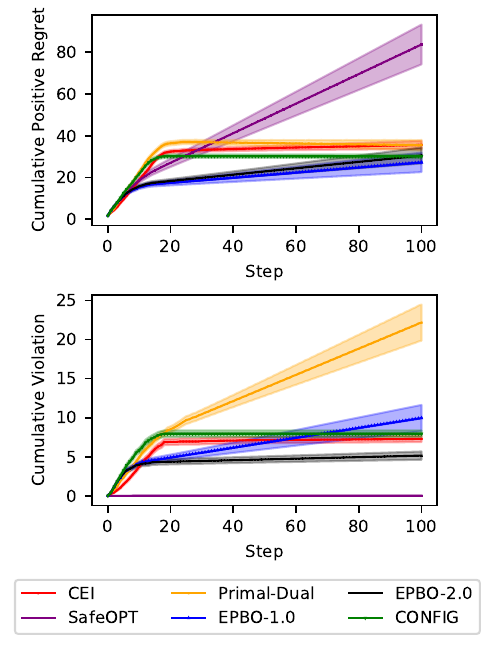}
\caption{\revReb{Cumulative regret and violation of different algorithms. The shaded area represents $\pm0.1\textsf{ standard deviation}$ and EPBO-$\rho$ represents EPBO with penalty $\rho$.}}
\label{fig:regret_vio}
\end{figure}
\begin{figure}[htbp!]
\centering
\includegraphics{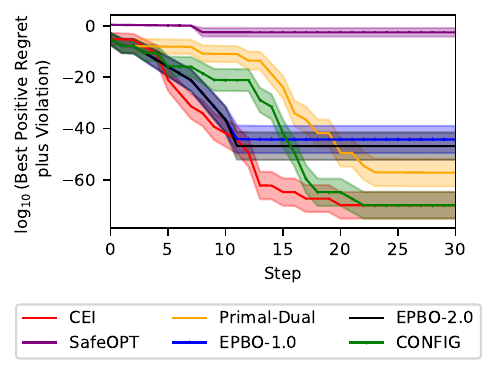}
\caption{\revReb{The best suboptimality plus violation up to step $t$, that is $\min_{\tau\in[t]}[f(x_\tau)-f^*]^++[g(x_\tau)]^+$, of different algorithms. The shaded area represents $\pm0.1\textsf{ standard deviation}$. Here, EPBO-$\rho$ represents the EPBO algorithm with penalty term $\rho$. The \textsf{SafeOPT} has a higher initial value due to the restriction of starting from a feasible solution.}}
\label{fig:subopt_p_v}
\end{figure}
As shown in Fig.~\ref{fig:regret_vio}, we observe that SafeOPT has almost linear cumulative regret, although achieving close-to-zero violation. This is due to the possibility of getting stuck in a local minimum for SafeOPT. In contrast, the primal-dual method can achieve competitive performance with the state-of-the-art constrained expected improvement~(CEI) method in terms of cumulative regret but suffers from almost linear cumulative violation. The EPBO method may also suffer from almost linear cumulative violation if the penalty term is not set large enough. Our \textsf{CONFIG} algorithm can achieve performance competitive with the state-of-the-art constrained expected improvement method in terms of both cumulative regret and cumulative violation. 


We further show the evolution of the best-observed suboptimality plus violation up to the current step $t$ in Fig.~\ref{fig:subopt_p_v}, which measures the empirical convergence speed to the optimal feasible solution. We notice that our algorithm again achieves empirical convergence speed competitive with the state-of-the-art CEI method. 

We then test our algorithms on a set of infeasible instances, where we shift the sampled constraint function $g(x)$ by $\epsilon-\min_{x\in\mathcal{X}}g(x)$ with $\epsilon=0.1$. Over 50 such infeasible instances, our algorithm consistently declares infeasibility, within 16.3 steps on average.

\subsection{Artificial Numerical Instances}
\label{subsec:art_nume}
In this section, we present the results of a set of global optimization test problems. 
We consider the global optimization problem in two-dimensional space with $\mathcal{X}=[-10, 10]^2$.
We then introduce a set of potentially non-convex and multi-modal functions as in Tab.~\ref{tab:arti_funcs} in the Appendix~\ref{app:exp_funcs}. Using these functions as the objective or constraint, we construct a set of global optimization benchmark problems as shown in Tab.~\ref{tab:arti_probs}. We use the common squared exponential kernel for all the experiments.  

We use constrained regret defined as
\bee
\min_{\tau\in[t]}[f(x_\tau)-f^*]^++[g(x_\tau)]^+, 
\eee
to measure the quality of the solution found by different algorithms up to step $t$. It can be seen that the smaller the constrained regret is, the closer to the constrained optimum our solution is in terms of optimality and violations.  
\begin{figure*}[htbp!]
\centering
\includegraphics[width=\textwidth]{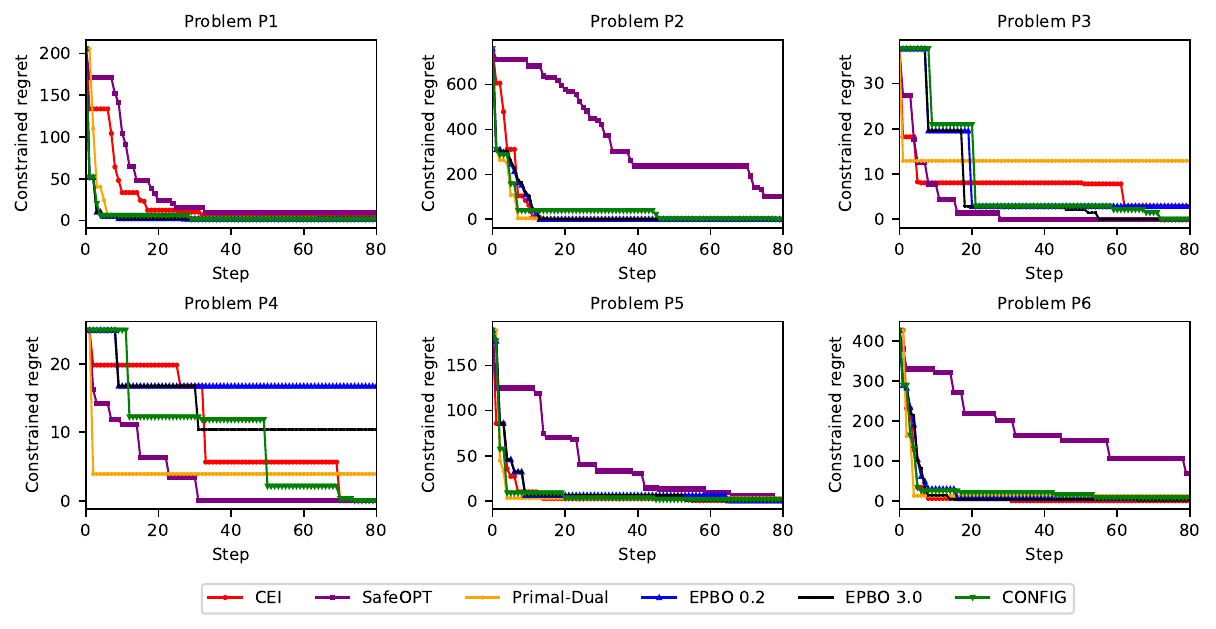}
\vspace{-0.3cm}
\caption{Convergence of constrained regret for the collections of problems in Tab.~\ref{tab:arti_probs}.}
\label{fig:P_collect}
\end{figure*}
\begin{table}[htbp!]
\centering
\caption{Artificial numerical problems constructed using the functions from Tab.~\ref{tab:arti_funcs}, where $\textsf{Qr}_\mathcal{X}(h)\defeq\frac{3}{4}\min_{x\in\mathcal{X}}h(x)+\frac{1}{4}\max_{x\in\mathcal{X}}h(x)$.}
\label{tab:arti_probs}
\begin{tabular}{lcc}
\hline $\text{ Problem }$ & \multicolumn{1}{c}{$\text{ Objective}$ $f$} & \multicolumn{1}{c}{$\text{ Constraint}$ $g$ } \\
\hline $\text{ P1 }$ & $\operatorname{Br}$ & $\operatorname{SinQ}-\textsf{Qr}_\mathcal{X}(\operatorname{SinQ})$  \\
$\text{ P2 }$ & $\operatorname{MBr}$ & $\operatorname{SinQ}-\textsf{Qr}_\mathcal{X}(\operatorname{SinQ})$  \\
$\text{ P3 }$ & $\operatorname{Br}$ & $\operatorname{InvBowl}-\textsf{Qr}_\mathcal{X}(\operatorname{InvBowl})$  \\
$\text{ P4 }$ & $\operatorname{MBr}$ & $\operatorname{InvBowl}-\textsf{Qr}_\mathcal{X}(\operatorname{InvBowl})$  \\
$\text{ P5 }$ & $\operatorname{Br}$ & $\operatorname{Bowl}-\textsf{Qr}_\mathcal{X}(\operatorname{Bowl})$  \\
$\text{ P6 }$ & $\operatorname{MBr}$ & $\operatorname{Bowl}-\textsf{Qr}_\mathcal{X}(\operatorname{Bowl})$  \\
\hline 
\end{tabular}
\end{table}

Fig.~\ref{fig:P_collect} demonstrates the convergence of constrained regret for different algorithms. The results show that the \textsf{CONFIG} algorithm consistently achieves a superior or comparable speed of finding the constrained optimal solution, as compared to the other state-of-the-art baselines. \revReb{We can also notice that in problem 3 and problem 4, SafeOPT outperforms both CEI and our method. Intuitively, it happens when the constrained global optimum lies in the same local feasible set as the initial points. In such a case, our method may spend lots of samples exploring globally, while SafeOPT identifies the local (and global) optimal solution quickly by restricting to the local feasible set.}   

\subsection{Tuning The P Controller of a Building}
\label{subsec:tune_P}
In this section, we consider a building temperature controller tuning problem. We use \textsf{Energym}~\cite{scharnhorst2021energym} with a lower level \textsf{Modelica}~\cite{fritzson1998modelica} model with a single thermal zone as our simulator. We consider a P controller controlling the temperature. The mathematical expression of the P controller is given as
$$
u=K_p(\textsf{set-point}-\textsf{temperature}),
$$
where $K_p$ is the proportional gain and $u$ is the heating control. We use $\theta$ to denote the tuning parameters, which include the proportional gain, daytime set-point for temperature, and the switching time from nighttime set-point to daytime set-point. Here, the nighttime set-point is fixed. Our objective is to minimize energy consumption while keeping the average temperature deviation below a given threshold. Our problem is formulated as in Eq.~\eqref{eqn:build_formulation}.
\bee\label{eqn:build_formulation}
\min_{\theta\in\Theta} \quad J(\theta) \quad\text{subject to}\quad   g(\theta)-g_{\mathrm{thr}}\leq 0,\quad
\eee
where $J(\theta)$ represents energy consumption in this example, $g(\theta)$ represents average temperature deviation and $g_\mathrm{thr}$ is a given average temperature deviation threshold. 
In each step, our algorithm selects a set of promising controller parameters. We use this set of parameters to run a closed-loop building simulation over a period of $1$ day. 


Based on the closed-loop trajectory, we calculate the key performance indicators, that is, energy consumption and average temperature deviation. We then update the Gaussian process regression with the new data added.  
Since energy consumption and temperature deviations are in two different scales, we use normalized energy consumption $\frac{J}{\sigma_J}$ plus normalized average temperature deviation $\frac{[g-g_\mathrm{thr}]^+}{\sigma_g}$ to measure the quality of the solution, where $\sigma_J, \sigma_g$ are standard deviations in a set of randomly sampled data. 
\begin{figure}[h!]
\centering
\includegraphics{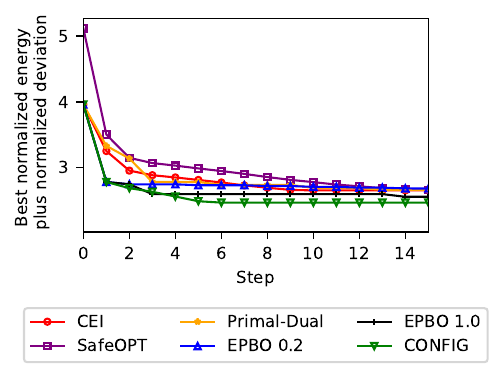}
\caption{Best normalized energy consumption plus normalized temperature deviation up to the current step. For \textsf{SafeOPT}, the initial required feasible solution is obtained based on domain knowledge.}
\label{fig:best_normalized_energy_plus_dev}
\end{figure}

Fig.~\ref{fig:best_normalized_energy_plus_dev} gives the best-normalized energy consumption plus average temperature deviation up to the current step. We observe that within $15$ steps, our algorithm identifies the solution with minimized normalized energy consumption plus normalized average temperature deviation, verifying the competitiveness of our algorithm with the currently popular constrained expected improvement method again. In contrast, the SafeOPT method can be overly cautious and fails to find a solution with performance comparable to ours.

\section{Conclusion}
In this paper, we have presented \algname{\OurAlg}, a simple and effective algorithm, for constrained efficient global optimization of expensive black-box functions by exploiting the principle of \emph{optimism in the face of uncertainty}. Specifically, our algorithm solves an auxiliary constrained optimization problem with the lower confidence bound~(LCB) surrogates as objective and constraints to get the next sample. Theoretically, we show that our algorithm enjoys the same cumulative regret bound as that in the unconstrained case. In addition, we show that the cumulative constraint violations have upper bounds in terms of maximum information gain, which are similar to the cumulative regret bounds of the objective function. For commonly used M\'atern~(\revReb{$\nu>\frac{d}{2}$}) and Squared Exponential kernels, our bounds are sublinear and allow us to derive a convergence rate to the optimal solution of the original constrained problem. Furthermore, our method naturally provides a scheme to declare infeasibility when the original black-box optimization problem is infeasible. Numerical experiments over sampled instances, artificial problems, and a building controller tuning problem, all corroborate the effectiveness of our algorithm, which is competitive with the popular CEI method. \revReb{As compared to existing state-of-the-art constrained efficient global/Bayesian optimization methods, \OurAlg{ significantly} improves the theoretical guarantees, while achieving competitive empirical performance.} 

One limitation of our work is that the auxiliary problem may be infeasible with misspecified hyperparameters, which needs to be addressed in future work.

\newpage
\bibliographystyle{icml2023} 
\bibliography{refs.bib}

\appendix
\section{Examples of kernel functions}
\label{app:ex_kerfuncs}
\begin{itemize}
\item\emph{Linear}: $$k(x,y)=x^Ty.$$
\item\emph{Squared Exponential}:  $$k(x,y) = \exp{\left\{-\frac{\norm{x-y}^2}{\ell^2}\right\}}.$$
\item\emph{M\'atern}: 
$$
\begin{aligned}
k&(x,y)=\\
&\frac{2^{1-\nu}}{\Gamma(\nu)}\left(\sqrt{2 \nu} \frac{\norm{x-y}}{\rho}\right)^{\nu} K_{\nu}\left(\sqrt{2 \nu} \frac{\norm{x-y}}{\rho}\right),    
\end{aligned}
$$ 
where $\rho$ and $\nu$ are positive parameters of the kernel function, $\Gamma$ is the gamma function, and $K_{\nu}$ is the modified Bessel function of the second kind.
\end{itemize}

\section{Proof of Corollary~\ref{cor:hpb}}
\begin{proof}
\rev{Note that $l_{i,t}(x)$ and $u_{i,t}(x)$ are random variables since they depend on the observations $y_{i,t}$'s, which are corrupted by random noise. The following events definitions are with respect to the randomness from noise. For simplicity of notation, we use $\{l_{i,t}(x)\leq y \leq u_{i,t}(x)\}$ to denote the intersection event that the random variable $l_{i,t}(x)\leq y$ and $u_{i,t}(x)\geq y$, where $y\in\mathbb{R}$, $i\in\{0\}\cup[N]$ and $t\geq1$. We define the event $\mathcal{E}_{i}\triangleq\cap_{x\in\mathcal{X}}\cap_{T\geq t\geq1}\{l_{i,t}(x)\leq g_i(x)\leq u_{i,t}(x)\}, \forall i\in[N]$ 
and the event $\mathcal{E}_{0}\triangleq\cap_{x\in\mathcal{X}}\cap_{T\geq t\geq1}\{l_{0,t}(x)\leq f(x)\leq u_{0,t}(x)\}$. We have,
\bse
\begin{align}
\mathbb{P}\left(\cap_{i=0}^N\mathcal{E}_i\right)&=1-\mathbb{P}\left(\overline{\cap_{i=0}^N\mathcal{E}_i}\right)\\
&=1-\mathbb{P}\left({\cup_{i=0}^N\overline{\mathcal{E}_i}}\right)\\
&\geq 1-\sum_{i=0}^N\mathbb{P}\left({\overline{\mathcal{E}_i}}\right)\label{eqn_line:union_bound}\\
&\geq 1-\sum_{i=0}^N\frac{\delta}{N+1}\label{eqn_line:lem}\\
&= 1-\delta,
\end{align}
\ese
where the inequality~\eqref{eqn_line:union_bound} follows by union bound of probability and the inequality~\eqref{eqn_line:lem} follows by Lemma~\ref{lem:conf_int}.
The conclusion then follows.
}

\end{proof}

\section{Proof of Lemma~\ref{lem:bound_ind_rv}}
\begin{proof}
By Corollary~\ref{cor:hpb}, with probability at least $1-\delta$, we have, for all $x\in\mathcal{X}$ and $1\leq t\leq T$,
\bse
\begin{eqnarray}
f(x)&&\in[l_{0,t}(x), u_{0,t}(x)]\enspace,\label{int:obj}\\
g_i(x)&&\in[l_{i,t}(x), u_{i,t}(x)]\enspace,\quad \forall i\in[N].\label{int:const}
\end{eqnarray}
\ese
All the following statements are conditioned on the above joint event happening, which has a probability of at least $1-\delta$.

\rev{By setting $x=x^*$ in~\eqref{int:const} and the feasibility of the optimal solution $x^*$, we have 
\bee
l_{i,t}(x^*)\leq g_i(x^*)\leq0, \quad\forall i\in[N].
\eee
} 

Therefore, infeasibility is not declared in line~\ref{alg_line:declare_inf} of the Alg.~\ref{alg:lcb2}, and $x^*$ is a feasible solution for the auxiliary problem in line~\ref{alg_line:aux_prob} of the Alg.~\ref{alg:lcb2}.
Therefore, we have
\bee
l_{0,t}(x_t)\leq l_{0,t}(x^*).~\label{ine:optim_imp}
\eee
Therefore, 
\bse
\begin{align}
r_t&=f(x_t)-f^*\\
&\leq r_t^+\\
&=[f(x_t)-l_{0,t}(x_t)+l_{0,t}(x_t)-f^*]^+\\
&\leq[f(x_t)-l_{0,t}(x_t)]^++[l_{0,t}(x_t)-f^*]^+\label{ine:pos_tri_ineq}\\
&\leq [u_{0,t}(x_t)-l_{0,t}(x_t)]^++[l_{0,t}(x_t)-l_{0,t}(x^*)]^+\label{ine:int_change_obj}\\
&= u_{0,t}(x_t)-l_{0,t}(x_t)\label{ine:opt_cond}\\
&= 2\betac_{0,t}\sigma_{0,t-1}(x_t),
\end{align}
\ese
where the inequality~\eqref{ine:int_change_obj} follows by the inequalities~\eqref{int:obj}, the inequality~\eqref{ine:pos_tri_ineq} follows by the fact that $[a+b]^+\leq[a]^++[b]^+,\forall a,b\in\mathbb{R}$, \rev{and the equality~\eqref{ine:opt_cond} follows by the inequality~\eqref{ine:optim_imp}}. 
We now consider the constraint violation,
\bse
\begin{align}
v_{i,t}&=[g_i(x_t)]^+\\
&=[g_i(x_t)-l_{i,t}(x_t)+l_{i,t}(x_t)]^+\\
&\leq [g_i(x_t)-l_{i,t}(x_t)]^++[l_{i,t}(x_t)]^+\\
&= [g_i(x_t)-l_{i,t}(x_t)]^+\label{eq:opt_cond}\\
&\leq[u_{i,t}(x_t)-l_{i,t}(x_t)]^+\label{ine:int_change}\\
&=2\betac_{i,t}\sigma_{i,t-1}(x_t),
\end{align}
\ese
where the equality~\eqref{eq:opt_cond} follows by the feasibility of $x_t$ for the auxilliary problem in line~\ref{alg_line:aux_prob} of Alg.~\ref{alg:lcb2}. 
\end{proof}

\section{Proof of Theorem~\ref{thm:R_bound}}
\begin{proof}
By the monotonicity of $\betac_{i,t}$ and Lem.~\ref{lem:bound_cumu_sd}, we have,
\bse
\begin{align}
R_T&= \sum_{t=1}^Tr_t\\
&\leq \sum_{t=1}^Tr_t^+=R_T^+\\
&\leq \sum_{t=1}^T2\betac_{0,t}\sigma_{0,t}(x_t)\\
&\leq 2\betac_{0,T}\sum_{t=1}^T\sigma_{0,t}(x_t)\\
&\leq 4\betac_{0,T}\sqrt{(T+2)\gamma_{0,T}},~\label{inq:bound_R_proof}
\end{align}
\ese
and $\forall i\in[N]$,
\bse
\begin{align}
\mathcal{V}_{i,T}&= \sum_{t=1}^Tv_{i,t}\\
&\leq \sum_{t=1}^T2\betac_{i,t}\sigma_{i,t}(x_t)\\
&\leq 2\betac_{i,T}\sum_{t=1}^T\sigma_{i,t}(x_t)\\
&\leq 4\betac_{i,T}\sqrt{(T+2)\gamma_{i,T}},~\label{inq:bound_V_proof}
\end{align}
\ese
Combining with the expression of $\beta_{i,t}$ in Definition~\ref{def:hpb_int}, we can derive
\begin{align*}
R_T&\leq R_T^+\leq 4\betac_{0,T}\sqrt{(T+2)\gamma_{0,T}}=\mathcal{O}(\gamma_{0,T}\sqrt{T})\enspace,\\
\mathcal{V}_{i,T}&\leq 4\betac_{i,T}\sqrt{(T+2)\gamma_{i,T}}=\mathcal{O}(\gamma_{i,T}\sqrt{T})\enspace,\quad\forall i\in[N],
\end{align*}
which completes the proof.
\end{proof}

\section{Proof of Theorem~\ref{thm:converge_opt}}
\begin{proof}
Consider the sum,
\begin{align*}
& \sum_{t=1}^T \left([f({x}_t)-f^*]^++\sum_{i=1}^N[g_i(x_t)]^+\right)\\
=& R_T^++\sum_{i=1}^N\mathcal{V}_{i,T}\\
\leq&4\sqrt{T+2}\sum_{i=0}^N\betac_{i,T}\sqrt{\gamma_{i,T}},
\end{align*}
where the last inequality follows by the inequalities~\eqref{inq:bound_R_proof} and~\eqref{inq:bound_V_proof}.  
Therefore,
\begin{align*}
&\min_{t\in[T]} \left([f({x}_t)-f^*]^++\sum_{i=1}^N[g_i(x_t)]^+\right) \\
\leq & \frac{1}{T}\sum_{t=1}^T \left([f({x}_t)-f^*]^++\sum_{i=1}^N[g_i(x_t)]^+\right)\\
\leq&\frac{4\sqrt{T+2}\sum_{i=0}^N\betac_{i,T}\sqrt{\gamma_{i,T}}}{T}.
\end{align*}
So there exists $\tilde{x}_T\in\{x_1, x_2, \cdots, x_T\}$, such that, 
$$
[f(\tilde{x}_T)-f^*]^++\sum_{i=1}^N[g_i(\tilde{x}_T)]^+\leq\frac{4\sqrt{T+2}\sum_{i=0}^N\betac_{i,T}\sqrt{\gamma_{i,T}}}{T}.
$$
Since $[f(\tilde{x}_T)-f^*]^+$ and $[g_i(\tilde{x}_T)]^+$ are non-negative, we have
\begin{align*}
f(\tilde{x}_T)-f^* &\leq [f(\tilde{x}_T)-f^*]^+\leq \frac{4\sqrt{T+2}\sum_{i=0}^N\betac_{i,T}\sqrt{\gamma_{i,T}}}{T}\\
[g_j(\tilde{x}_T)]^+&\leq \frac{4\sqrt{T+2}\sum_{i=0}^N\betac_{i,T}\sqrt{\gamma_{i,T}}}{T}\enspace,\quad\forall j\in[N].
\end{align*}
The desired conclusion then follows.
\end{proof}
\section{Proof of Theorem~\ref{thm:declare_inf}}
\begin{proof}


Suppose up to step $T$, infeasibility has not been declared. We have,
$$l_{i,t}(x_t)\leq0,\forall i\in[N],$$
and since $\min_{x\in\mathcal{X}}\max_{i\in[N]}g_i(x)=\epsilon$, there exists $i_t\in\arg\max_{i\in[N]}g_i(x_t)$, such that  
$$u_{i_t,t}(x_t)\geq g_{i_t}(x_t)\geq\epsilon$$
with probability at least $1-\delta$.

Accordingly, we have 
$$\sum_{t=1}^T\sum_{i\in[N]}2\beta_{i,t}\sigma_{i,t-1}(x_t)=\sum_{t=1}^T\sum_{i\in[N]}(u_{i,t}(x_t)-l_{i,t}(x_t))\geq \sum_{t=1}^T(u_{i_t,t}(x_t)-l_{i_t,t}(x_t))\geq\sum_{t=1}^T\epsilon=T\epsilon.$$
Meanwhile,
$$
\begin{aligned}
2\sum_{i\in[N]}\sum_{t=1}^T\beta_{i,t}\sigma_{i,t-1}(x_t)&\leq 2\sum_{i\in[N]}\beta_{i,T}\sum_{t=1}^T\sigma_{i,t-1}(x_t)\\
&\leq 2\sum_{i\in[N]}\beta_{i,T}\sqrt{4(T+2) \gamma_{i,T}},
\end{aligned}
$$
where the last inequality follows by Lem.~\ref{lem:bound_cumu_sd}. Therefore,
$$2\sum_{i\in[N]}\beta_{i,T}\sqrt{4(T+2) \gamma_{i,T}}\geq T\epsilon,$$  
which implies $\epsilon\leq\frac{2\sum_{i\in[N]}\beta_{i,T}\sqrt{4(T+2) \gamma_{i,T}}}{T}$. 
By Def.~\ref{def:hpb_int}, we have 
$$\frac{2\sum_{i\in[N]}\beta_{i,T}\sqrt{4(T+2) \gamma_{i,T}}}{T}=\mathcal{O}\left(\sum_{i\in[N]}\frac{\gamma_{i,T}}{\sqrt{T}}\right).$$
Hence, $\exists \tilde{C}>0$, such that,
\bee\epsilon\leq \tilde{C}\sum_{i\in[N]}\frac{\gamma_{i,T}}{\sqrt{T}}.\eee
That is 
\bee C\epsilon\leq\sum_{i\in[N]} \frac{\gamma_{i,T}}{\sqrt{T}},\label{inq:nece_cond_not_declare}\eee
where $C=\frac{1}{\tilde{C}}$.
However, since $\lim_{T\to\infty}\frac{\gamma_{i,T}}{\sqrt{T}}=0,\forall i\in [N]$, the inequality~\eqref{inq:nece_cond_not_declare} will be violated when $T$ is large enough. So infeasibility will be declared on or before the first time the above inequality is violated, which is $\overline{T}$.      
\end{proof}
\section{Explicit functions for the artificial numerical instances}
\label{app:exp_funcs}
\begin{table*}[htbp]
\renewcommand{\arraystretch}{1.5}
\centering
\begin{adjustbox}{angle=0}
\begin{tabular}{l|l}
\hline  { Name } & \multicolumn{1}{c}{{ Function }} \\
\hline { Branin (Br) } & $\operatorname{Br}\left(x_1, x_2\right)=\left(x_2-\frac{5.1}{4 \pi^2} x_1^2+\frac{5}{\pi} x_1-6\right)^2+10\left(1-\frac{1}{8 \pi}\right) \cos \left(x_1\right)+10$ .  \\
\hline $\begin{array}{l}
\text{ Modified } \\
\text{ Branin (MBr) }
\end{array}$ & $\operatorname{MBr}\left(x_1, x_2\right)=\operatorname{Br}\left(x_1, x_2\right)+20 x_1-30 x_2$ .  \\
\hline { Bowl (Bowl) } & 
$\begin{array}{l}
\text{ Bowl }\left(x_1, x_2\right)=\frac{1}{2}\left(\left\|x-c_{\text{bowl }}\right\|^2-R_{\text{bowl }}^2\right) \text{ where } R_{\text{bowl }}=10 \text{ and } \\
c_{\text{bowl }}=[-3,-3] .
\end{array}$ \\
\hline $\begin{array}{l}
\text{ Inverted bowl } \\
\text{ (InvBowl) }
\end{array}$ & $\operatorname{InvBowl}\left(x_1, x_2\right)=-\operatorname{Bowl}\left(x_1, x_2\right)$ .  \\
\hline $\begin{array}{l}
\text{ Sine- } \\
\text{ quadratic } \\
\text{ (SinQ) }
\end{array}$ & $\operatorname{SinQ}\left(x_1, x_2\right)=\sin \left(\frac{x_1^2+x_2^2}{10}\right)$ .  \\
\hline
\end{tabular}
\end{adjustbox}
\caption{Explicit functions used to construct the artificial numerical instances.}
\label{tab:arti_funcs}
\end{table*}

\end{document}